\documentclass[11pt]{article}
\usepackage{amssymb}
\usepackage{amscd}
\usepackage{epsf}
\newtheorem{thm}{Theorem}[section]
\newtheorem{lem}[thm]{Lemma}

\newtheorem{defi}[thm]{Definition}

\newcommand{\inv}{^{-1}}

\newcommand{\C}{{\mathbb{C}}}
\newcommand{\Z}{{\mathbb{Z}}}
\newcommand{\R}{{\mathbb{R}}}
\newcommand{\N}{{\mathbb{N}}}
\newcommand{\T}{{\mathbb{T}}}

\newcommand{\Gm}{{\mathbb{G}}_m}

\newcommand{\Hom}{\mbox{\rm Hom}}
\newcommand{\diag}{\mbox{\rm diag}}
\newcommand{\trop}{\mbox{\rm trop}}

\begin{document}

\title{A tropical view on the Bruhat-Tits building of SL and its compactifications}
\author{Annette Werner\\
\small Institut f\"ur Mathematik\\[-1.5ex]
\small Goethe-Universit\"at Frankfurt\\[-1.5ex]
\small Robert-Mayer-Strasse 8\\[-1.5ex]
\small  D- 60325 Frankfurt\\[-1.5ex]
\small email: werner@math.uni-frankfurt.de}

\date{}

\maketitle

\centerline{\bf Abstract:} We describe the stabilizers of points in the Bruhat-Tits building of the group $SL$ with tropical geometry. There are several compactifications of this building associated to algebraic representations of $SL$. We show that the fans used to compactify apartments in this theory are given by tropical Schur polynomials.
\small 
\\[0.3cm]

\centerline{{\bf 2000 MSC: }20E42, 20G25, 52B99 } 
\centerline{{\bf Keywords: }Bruhat-Tits buildings, tropical hypersurfaces, Schur polynomials}

\normalsize 
\section*{Introduction}
In this paper we relate  tropical geometry to point stabilizers and compactifications of the  Bruhat-Tits building associated to the group $SL$.

Let $V$ be  vector space of dimension $n$ over a field $K$ with a non-trivial non-Archimedean absolute value. Then the Bruhat-Tits building $\mathfrak{B}(SL(V), K)$ is a metric space with a continuous $SL(V)(K)$-action. It consists of apartments, i.e. real cocharacter spaces of maximal split tori. They are glued together with the help of certain subgroups of $SL(V)(K)$ which turn out to be the stabilizers of points in the building. 
We can regard the apartments which are real vector spaces from a tropical point of view. Every element in the special linear group over $K$ defines a real matrix by taking the negative valuation of all coefficients. Replacing addition and scalar multiplication by tropical addition and tropical scalar multiplication, such a matrix defines a tropical linear map. In theorem 1.2 we prove that the tropical linear stabilizer of a point in an apartment coincides with its stabilizer with respect to the $SL(V)(K)$-action on the  Bruhat-Tits building. 

Every geometrically irreducible algebraic representation of $SL(V)$ induces a compactification of the building $\mathfrak{B}(SL(V),K)$ by \cite{we3}. This compactification is defined by glueing together compactified apartments. Every apartment is compactified with the help of a fan defined by the combinatorics of the weights given by the representation. We assume that $K$ has characteristic zero. Every representation of $SL(V)$ is then given by a partition $\lambda = (\lambda_1 \geq \ldots \geq \lambda_n \geq 0)$, and its character is the Schur polynomial $S_\lambda$ associated to this partition. We look at the tropical hypersurface defined by the Schur polynomial which gives rise to a natural polyhedral complex. In Theorem 2.4 we prove that the following fact. If the coefficients of the Schur polynomial all have valuation zero, then the fan used to compactify one apartment is precisely the fan given by the tropical Schur polynomial. 

This result was motivated by the following example. Consider the identical representation of $SL(V)$. The corresponding partition is $(1, 0 \ldots, 0)$, which leads to the Schur polynomial $S_{(1,0,\ldots, 0) }(z_1, \ldots, z_n)= z_1 + \ldots + z_n$.  In this case the fan used to compactify
the apartment plays an important role in tropical convexity defined, as defined by Develin and Sturmfels in \cite{dest}. 

This paper was inspired by \cite{jsy} who took another look the Bruhat-Tits building for $SL$ from a tropical point of view by describing and computing convex subsets of the building via tropical geometry. 

{\bf Acknowledgements: }I thank Michael Joswig and Thorsten Theobald for useful discussions on tropical geometry. 

\section{Stabilizers}

\subsection{The tropical torus} Tropical geometry is based on the tropical semiring $(\mathbb{R}, \oplus, \odot)$ with  $a \oplus b = \max\{a,b\}$ and $a \odot b = a +b$. Some authors use  the $(\min, +)$-version of the tropical semiring instead of the $(\max,+)$ version.

The space $\R^n$ together with componentwise addition $\oplus$ is a semimodule under the semiring $(\R, \oplus, \odot)$, if we put $a \odot (x_1,\ldots, x_n) = (a+ x_1, \ldots, a+ x_n)$ for $a \in \R$ and $(x_1,\ldots, x_n) \in \R^n$.  Let
\[ \T^{n-1} = \R^n / \R (1,\ldots, 1)\]
be the quotient
of $\R^n$ after the following equivalence relation:  $(x_1,\ldots, x_n) \sim (y_1,\ldots, y_n)$ if and only if there exists some $a \in \R$ such that $x_i = a \odot y_i = a + y_i $ for all $i$. We call $\T^{n-1}$ the tropical torus of rank $n-1$ as in \cite{jo}. We endow $\T^{n-1}$ with the quotient topology. In the literature, the space $\T^{n-1}$ is often called tropical projective space. We avoid this terminology, since $\T^{n-1}$ is not compact. In section 2.2, we consider a compactification of $\T^{n-1}$, which could be viewed as a tropical analogue of the projective space. 

Let $K$ be a field with a valuation map $v: K^\times \rightarrow \R$. We put $v(0) = \infty$. 
Consider the  negative valuation  map $-v: K \rightarrow \mathbb{R}_{- \infty} = \mathbb{R} \cup \{-\infty\}$. 
We extend the operations $\odot$ and $\oplus$ to $\R_{-\infty}$ by  $a \odot - \infty = - \infty$ for all $a \in \R_{-\infty}$ and $a \oplus -\infty = a$ for all $a \in \R_{-\infty}$. 

\begin{defi}
Let $g = (g_{ij})_{i,j}$ be a matrix in $GL_n(K)$. Then we define the associated tropical matrix by
$g_{\trop} = (-v ( g_{ij}))_{i,j} \in \mbox{Mat}_{n \times n} (\R)$. For every vector $x = ~^t(x_1,\ldots, x_n) \in \R^n$ we define the 
  vector $ g_{\trop} \cdot x = ~^t(y_1, \ldots, y_n) \in \R^n$ by the tropicalized linear action
 \[ y_i = -v(g_{i1}) \odot x_1 \oplus \ldots \oplus -v(g_{in}) \odot x_n = \max_j \{ -v(g_{ij}) + x_j\}.\]
\end{defi}

Note that $y_i$ lies indeed in $\R$, since at least one entry $g_{ij}$ in the $i$-th line must be non-zero, so that at least for one $j$ the term $-v(g_{ij}) + x_j$ is not equal to $-\infty$. 

Beware that this does not define an action of $GL_n(K)$ on $\R^n$, as the following example shows. Take $n = 2$, and put
\[g = \left( \begin{array}{ll} 1 & 1\\ 0 & 1 \end{array} \right) \mbox{ and } h = \left( \begin{array}{ll} 1 & 0 \\ -1 & 1 \end{array} \right). \mbox{ Then}\]
\[(gh)_{\trop} \cdot \left(\begin{array}{l} x_1 \\ x_2 \end{array} \right) = \left(\begin{array}{l} x_2 \\ \max\{x_1, x_2\} \end{array} \right), \mbox{ but}\]
\[g_{\trop} \cdot \left( h_{\trop} \cdot \left(\begin{array}{l} x_1 \\ x_2 \end{array} \right)\right) = \left(\begin{array}{l} \max\{x_1, x_2\} \\ \max\{x_1, x_2\}\end{array} \right).\]

\subsection{Bruhat-Tits buildings} Now we want to explain how the tropical torus $\T^{n-1}$  appears in the context of Bruhat-Tits buildings. Let us first introduce some notation.  We fix a ground field $K$  with a non-trivial non-Archimedean absolute value $| \quad |$, i.e. there is a valuation map $v : K^\times \rightarrow \R$ such that $|x| = e^{-v(x)}$ for all $x \in K^\times$. For example, $K$ could be a local field, i.e. a field which is complete with respect to discrete valuation and has finite residue field. In this case, $K$ is either a finite extension of $\mathbb{Q}_p$ or a field of formal Laurent series over a finite field. Alternatively, $K$ could be $\C_p$, which denotes the completion of the algebraic closure of $\mathbb{Q}_p$,  or the field of Puiseux series over $\C$ which is quite popular in tropical geometry. 
By $\mathcal{O}_K$ we denote the ring of integers in $K$, i.e. $\mathcal{O}_K = \{ x \in K: |x| \leq 1\}$. 

Let $G$ be a reductive group over $K$. In the two groundbreaking papers \cite{brti1} and \cite{brti2} Bruhat and Tits  define a metric space $\mathcal{B}(G,K)$, now called Bruhat-Tits building, endowed with a continuous $G(K)$-action, if $G$ and $K$ satisfy certain assumptions. These  assumptions are for example fulfilled if the ground field is a local field or if the group $G$ is split, i.e. if $G$ contains a maximal torus which is split over $K$.  If the valuation is discrete, the building $\mathfrak{B}(G,K)$ carries  a polysimplicial structure.  If the valuation is not discrete, there is still a notion of facets, see \cite{rou}. 

Let us now recall the construction of $\mathfrak{B}(G,K)$  in the special case $G = SL(V)$ where  $V$ is a vector space of dimension $n$ over $K$. 
Let $T$ be a maximal split torus in $SL(V)$, i.e. there is a basis $e_1,\ldots, e_n$ of $V$ such that $T$ consists of all diagonal endomorphisms in $SL(V)$ with respect to this basis.  Let $X^\ast(T) = \Hom_{K} (T, \Gm)$ denote the character group. Here we write $\Hom_K$ for morphisms of affine $K$-group schemes. Generally, we regard the algebraic groups $T$, $SL(V)$ etc. as group schemes over $K$ and write $T(K)$, $SL(V)(K)$ etc. for the groups of $K$-rational points (i.e. for the corresponding groups of matrices with coefficients in $K$). 

 We define $a_i \in X^\ast(T)$ by
\[a_i (\diag(s_1, \ldots, s_n)) = s_i,\]
where $\diag(s_1, \ldots, s_n)$ is the diagonal matrix with entries $s_1, \ldots, s_n$. 
Then $X^\ast(T) = \bigoplus_{i=1}^n \Z a_i / \Z (a_1 + \ldots + a_n)$. 

Besides, let $X_\ast(T) = \Hom_{K} (\Gm, T)$ denote the cocharacter group of $T$. 
Let $\eta_i : \Gm \rightarrow GL(V)$ be the cocharacter of $GL(V)$ mapping $x$ to the diagonal matrix with entries $t_1,\ldots, t_n$, where $t_i= x$ and $t_j = 1$ for $j \neq i$. Then 
\[X_\ast(T) = \{ m_1 \eta_1 + \ldots + m_n \eta_n : m_i \in \Z \mbox{ with }\sum_i m_i = 0\}.\] 
There is a natural pairing $X^\ast(T) \times X_\ast (T) \rightarrow \Z$ induced by $a_i(\eta_j)= \delta_{ij}$. 

The $\R$-vector space $A = X_\ast(T) \otimes_\Z \R = \{\sum_{i=1}^n x_i \eta_i : x_i \in \R \mbox{ with } \sum_i x_i = 0\}$ is the apartment given by the torus $T$ in the Bruhat-Tits building $\mathcal{B}(SL(V), K)$. 
Mapping $\eta_1,\ldots, \eta_n$ to the canonical basis of $\R^n$ provides a homeomorphism
\[A \longrightarrow \T^{n-1}\]
between the apartment $A$ and the tropical space $\T^{n-1}$. In the course of this paper we always identify $A$ with $\T^{n-1}$ using the map above. 

For every  $t = \diag(t_1,\ldots, t_n) \in T(K)$  with entries $t_1,\ldots, t_n \in K^\times$ we define a point in $A$ by $\nu(t) = -v(t_1) \eta_1 + \ldots + -v(t_n) \eta_n$. Then $ t \in T(K)$ acts on $A$ by translation with $\nu(t)$. Besides, let $N$ be the normalizer of $T$ in $SL(V)$. For every element $n \in N(K)$  there is a  permutation $\sigma$ on $\{1, \ldots, n\}$ such that  $n(e_i) = t_{i} e_{\sigma(i)}$ for suitable $t_1, \ldots, t_n \in K^\times$. The Weyl group $W = N(K) / T(K)$ can therefore be identified with the 
symmetric group $\mathcal{S}_n$. Hence $W$ acts in a natural way on $A$ by permuting the coordinates of a given point. We can put both actions together to an action of $N(K)$ on $A$ by affine-linear transformations. If $A$ is endowed with its real topology and $N(K)$ is endowed with the topology
given by the topology on $K$, this action is continuous. 

Bruhat-Tits theory provides  a filtration on the $K$-rational points of the root groups for any reductive group.
In our special case this structure can be described as follows. The root system $\Phi(T, SL(V))$ consists of all characters of the form
$a_{ij} = a_i / a_j$ for all $i,j \in \{1,\ldots, n\}$ with $i \neq j$. The root group $U_a$ corresponding to the root $a = a_{ij}$ consists of all $u \in SL(V)$ such that 
\begin{eqnarray*}
u(e_k)& = &e_k \mbox{ for all }k \neq i \mbox{ and}\\
u(e_j) & = & e_j + \omega e_i
\end{eqnarray*}
for some $\omega \in K$. 

Then we have a homomorphism 
\[\psi_a: U_a(K) \rightarrow \Z \cup \{\infty\}\]
mapping $u$ to  $v(\omega)$. The Bruhat-Tits filtration on $U_a(K)$ is given by the subgroups
\[U_{a,l} = \{ u \in U_a(K): \psi_a(u) \geq l\}\]
for all $l \in \Z$. 

For all $x \in A$ consider the group  $U_x$ generated by $U_{a,-a(x)}$ for all $a \in \Phi(T,SL(V))$, and let $N_x$ be the subgroup of $N(K)$ consisting of all elements stabilizing $x$. Then the group $P_x = N_x U_x = U_x N_x$ can be used to define the Bruhat-Tits building $\mathfrak{B}(SL(V),K)$ as the quotient of $SL(V)(K) \times A$ after the following equivalence relation:
\[
(g,x) \sim (h,y),  \begin{array}{l}
\mbox{ if and only if  there exists an element } n \in N(K) \\
 \mbox{ such that } \nu(n)x= y \mbox{ and } g\inv h n \in P_x.
 \end{array}
 \]
The quotient space $\mathcal{B}(SL(V),K)$ is endowed with the product-quotient topology and admits a natural continuous $SL(V)(K)$-action via left multiplication in the first factor. For all $x \in A$, the group $P_x$ is the stabilizer of $x$ in $\mathcal{B}(SL(V),K)$. All subsets of $\mathcal{B}(SL(V),K)$ of the form $g A$ for $g \in SL(V)(K)$ are called apartments.

If the valuation on $K$ is discrete, then  $\mathcal{B}(PGL(V),K)$ can be identified with the Goldman-Iwahori space of all non-Archimedean norms on $K^n$ modulo scaling, see \cite{brti3}. Here a map $\gamma: K^n \rightarrow \mathbb{R}_{\geq 0}$ is a non-Archimedean norm if $\gamma(\lambda v) = |\lambda| \gamma (v)$ and $\gamma (v+w) \leq \sup\{\gamma(v),\gamma(w)\}$ for all $\lambda \in K$ and $v,w \in K^n$ and $\gamma(v) = 0$ only for $v=0$ holds. 
Via this identification, the apartment $A$ consists of all norms (modulo scaling) of the form
\[\gamma((\lambda_1, \ldots, \lambda_n))= \sup\{|\lambda_1| r_1, \ldots, |\lambda_n| r_n\} \]
for some real vector $(r_1,\ldots, r_n)$. 

In this case the simplicial complex $\mathfrak{B}(G,K)$ can be described as a flag complex whose vertex set consists of all  similarity classes of $\mathcal{O}_K$-lattices in $K^n$. Here two lattice classes are adjacent  if and only if there are representatives $M$ and $N$ of these two classes satisfying $\pi M \subset N \subset M$, where $\pi$ is a prime element in the ring of integers $\mathcal{O}_K$. 

\subsection{ A tropical view on the stabilizer $P_x$ } We will now show that the stabilizer groups $P_x$ defined above with Bruhat-Tits theory have a tropical interpretation.  Namely, we show that $P_x$ coincides with the set of all $g \in SL_n(K)$ stabilizing  $x$ under tropical matrix multiplication. Recall that the map $(g,x) \mapsto g_{\trop} \cdot x$ does not define an action of $SL_n(K)$ on $\R^n$. Hence the theorem below also shows that the tropical stabilizer of a point in $\R^n$ is a group which is not a priori clear. We use the basis $e_1, \ldots, e_n$ of $V$ to identify $SL(V)$ with $SL_n$.

\begin{thm} Let $x = \sum_i x_i \eta_i$ be a point in the apartment $A$, and let $P_x$ be the stabilizer of $x$ with respect to the action of $SL_n(K)$ on $\mathcal{B}(SL_n,K)$.  Then we have
\[P_x = \{ g \in SL_n(K): g_{\trop} \cdot ~^t(x_1, \ldots, x_n) = ~^t(x_1, \ldots, x_n)\}.\]
\end{thm}

{\bf Proof: }First we consider the case that every  $x_i$ lies in the image of the valuation map $v: K^\times \rightarrow \R$. Hence there exist elements $t_i \in K^\times $ satisfying $x_i = - v(t_i)$. For  $g= (g_{ij})_{i,j}$ in $SL_n(K)$ we have
\[g_{\trop} \cdot ~^t(x_1,\ldots, x_n) = ~^t( \max_j \{-v(g_{1j} t_j)\}, \ldots, \max_j \{-v(g_{nj} t_j)\}).\]
Hence $g_{\trop} ~^t(x_1,\ldots, x_n) = ~^t(x_1,\ldots, x_n)$ is equivalent to the fact that $\max_j \{-v( g_{ij} t_j  t_i^{-1})\} = 0$ for all  $i$.  Now $g_{ij}  t_i^{-1} t_j$ is the entry at position $(i,j)$ of the matrix $ t^{-1} g t$, where $t $ denotes the diagonal matrix with entries $t_1, \ldots, t_n$. Therefore the matrix $t^{-1} g t$ lies in $SL_n(\mathcal{O}_K)$. Conversely, every $g$ such that $t^{-1} g t$ lies in $SL_n({\mathcal{O}}_K)$ stabilizes $~^t(x_1,\ldots, x_n)$. Therefore 
\[\{g \in SL_n(K): g_{\trop} \cdot ~^t(x_1,\ldots, x_n) = ~^t(x_1,\ldots, x_n)\} = t SL_n(\mathcal{O}_K) t^{-1}.\]
 By Bruhat-Tits theory, $t SL_n({\mathcal{O}}_K) t^{-1}$ is indeed the stabilizer of $x = \sum_i -v(t_i) \eta_i$ in the apartment $A$. 

For a general point $x = \sum_i x_i \eta_i \in A$ there exists a non-Archimedean extension field $L$ of $K$ such that all $x_i$ are contained in the image of the valuation map of $L$ and a continuous $SL_n(K)$-equivariant embedding $\mathfrak{B}(SL_n, K) \hookrightarrow \mathfrak{B}(SL_n, L)$, see e.g. \cite{rtw1}, (1.2.1) and (1.3.4). We have just shown that $\{ g \in SL_n(L): g_{\trop} \cdot ~^t(x_1, \ldots, x_n) = ~^t(x_1, \ldots, x_n)\}$ is equal to the stabilizer of $x$ in the building over $L$. Intersecting with $SL_n(K)$ our claim follows.\hfill$\Box$

\section{Tropical Schur polynomials and compactifications of the building}

In this section we assume that our non-Archimedean field $K$ has characteristic zero. 

\subsection{Tropical Schur polynomials} 
Schur polynomials play an important role in the representation theory of the symmetric group and of the general linear group. The polynomial representations of $GL_n$ over $\C$ were determined by Isaai Schur in his thesis. See \cite{pro} for a modern account. In \cite{gr}, the theory is developped over arbitrary infinite ground fields. We recall the following facts. 

Fix a natural number $n \geq 2$. 
We call $\lambda = (\lambda_1, \ldots, \lambda_n)$ a partition (of $\lambda_1 + \ldots + \lambda_n$) if $\lambda_1 \geq \ldots \geq \lambda_n \geq 0$.

The Schur polynomial $S_\lambda$  can be defined by the following formula:
\[S_{\lambda}(z_1,\ldots, z_n) = \frac{ \det((z_j^{\lambda_i + n -i})_{i,j = 1, \ldots n})}{\det ((z_j^{n-i})_{i,j = 1, \ldots,n})}.\]
For properties of these homogeneous symmetric polynomials see \cite{mac}, I.3.

{\bf Example: }$S_{1,0, \ldots, 0}(z_1, \ldots, z_n) = z_1 + \ldots + z_n$.

Let $V$ be a vector space over $K$ of dimension $n$. For every partition $\lambda$ there is a geometrically
irreducible algebraic representation $\mathfrak{S}_\lambda(V)$ of $SL(V)$ such that for every $g \in SL(V)$ with eigenvalues $c_1,\ldots, c_n$ the trace of $g$ on $\mathfrak{S}_\lambda(V)$ is equal to $S_\lambda(c_1, \ldots, c_n)$, where $S_\lambda$ is the Schur polynomial associated to $\lambda$. 
The partitions $(\lambda_1, \ldots, \lambda_n)$ and $(\lambda_1+m, \ldots, \lambda_n+m)$ give the same representation for all integers $m$, and every geometrically irreducible algebraic representation of $SL(V)$ arises in this way. 

Let $W$ be a $K$-vector space and let $\rho: SL(V) \rightarrow GL(W)$ be an irreducible algebraic representation. Then the action of the torus $T$ on $W$ defines a decomposition of $W$ in weight spaces. Put
\[ W_\mu = \{ w \in W : \rho(t) (w) = \mu(t) w \mbox{ for all } t \in T(K)\}.\]
Then every character $\mu$ of $T$  with $W_\mu \neq 0$ is called a weight of $\rho$. 

A subset $\Delta$ of the root system $\Phi(T, SL(V))$ is called a basis if it is linear independent  and if every root $a$ in $\Phi(T, SL(V))$  can be written as $a = \sum_{c \in \Delta} m_c c$ where  either all $m_c$ are  non-negative or all $m_c$  are non-positive.  The Weyl group $W$ acts simply transitively on the set of bases, see \cite{bou}, chapter VI, \S1.5. 

For every $\Delta$ there is a highest weight $\mu_0(\Delta)$ such that for all weights $\mu$ we have  $\mu_0(\Delta) -\mu = \sum_{a \in \Delta} n_a a$ with non-negative integral coefficients $n_a$. 

Since the Schur polynomial $S_\lambda$  is the trace of the irreducible representation $\mathfrak{S}_\lambda(V)$, we can write it as 
\[S_\lambda(z_1,\ldots, z_n) = \sum_\mu \mbox{dim} (\mathfrak{S}_\lambda(V))_\mu z^\mu\]
with $z^\mu = z_1^{\mu_1} \cdots z_n^{\mu_n}$. 

Hence for any $\mu_1, \ldots, \mu_n \in \N_0$ with $\sum_i \mu_i = \sum_i \lambda_i$, the element $\mu_1 a_1  + \ldots + \mu_n a_n$ (modulo $\Z \sum_i a_i$) is a weight of the representation $\mathfrak{S}_\lambda(V)$ if and only if the monomial $z_1^{\mu_1} \cdots z_n^{\mu_n}$ occurs with a non-zero coefficient in the Schur polynomial $S_\lambda$.

The coefficients of the Schur polynomial $S_\lambda$ are given by Kostka numbers. To be precise, we have
\[S_\lambda = M_\lambda + \sum_{\mu} K_{\lambda \mu} M_\mu,\]
where $\mu$ runs over all partitions with $\sum_i \mu_i = \sum_i \lambda_i$, and where $M_\mu$ is the sum of the monomial $z_1^{\mu_1} z_2^{\mu_2} \cdots z_n^{\mu_n}$ and all the monomials we get from it by permuting the variables. The Kostka number $K_{\lambda \mu}$ denotes the number of ways one can fill the Young diagram associated to $\lambda$ with $\mu_1$ 1's, $\mu_2$ 2's etc. such that the entries are strictly increasing in the columns and non-decresasing in the rows.  Since $K_{\lambda \mu} \neq 0 $ implies that $\sum_{i=1}^k \lambda_i \geq \sum_{i=1}^k \mu_i$ for all $k$, it follows that the weight $\lambda_1 a_1 + \ldots +\lambda_n a_n$ is the highest weight of the representation $\mathfrak{S}_\lambda(V)$ with respect to the basis $\Delta= \{a_{12}, \ldots, a_{n-1 n}\}$. 

\begin{defi} Let $\lambda = (\lambda_1 \geq \ldots \geq \lambda_n \geq 0 )$ be a partition. By $\mathcal{T}(S_\lambda)$ we denote the tropical hypersurface of $\R^n$ given by the Schur polynomial
$S_\lambda$, i.e. $\mathcal{T}(S_\lambda)$ is the closure of the set
\[\{ (-v(\alpha_1), \ldots -v(\alpha_n)) \in \R^n: \alpha_1, \ldots, \alpha_n \in \overline{K}^\times \mbox{ with }S_\lambda(\alpha_1, \ldots, \alpha_n) = 0 \}. \]
\end{defi} 
Here we extend the valuation map to the algebraic closure $\overline{K}$ of $K$. Since $S_\lambda$ is homogeneous, ist also defines a tropical hypersurface in $\T^{n-1}$, which is the image of $\mathcal{T}(S_\lambda)$ under the quotient map $\R^n \rightarrow \T^{n-1}$. 

We write $S_\lambda (z_1, \ldots, z_n) = \sum_N c_N z^N$ with $N= (\nu_1, \ldots, \nu_n)$ and $z^N =z_1^{\nu_1} \cdots z_n^{\nu_n}$. (Here we write out all monomials.) Then by \cite{ekl}, theorem 2.1.1 we have
\begin{eqnarray*}
\lefteqn{\mathcal{T}(S_\lambda) = \{ (x_1, \ldots, x_n) \in \R^n:} \\
~ &  \max_{N= (\nu_1, \ldots, ´\nu_n)} \{ -v(c_N) + \nu_1 x_1 + \ldots + \nu_n x_n\} \mbox{ is attained at least twice}\}.
\end{eqnarray*}

The tropical hypersurface $\mathcal{T}(S_\lambda)$ gives rise to a polyhedral complex in $\R^n$. Let us recall the definition which works for arbitrary tropical hypersurfaces. 

For basics and terminology concerning polyhedral complexes see \cite{zi}. 
Consider the extended Newton polytope $P$ in $\R^{n+1}$ given by $S_\lambda = \sum_N c_N z^N$, i.e. $P$ is the convex hull of all points $(-v(c_N), \nu_1, \ldots, \nu_n)$ such that $c_N \neq 0$ in $\R^{n+1}$. The normal fan $\mathcal{N}_P$ associated to $P$ is the fan in $\R^{n+1}$ consisting of all cones 
\[\mathcal{N}_P(F) = \{w = (w_0, \ldots, w_n) \in  \R^{n+1}: F \subset \{x \in P:  <w,x> = \max_{y \in P} <w,y> \}  \},\]
where $F$ runs over the non-empty faces of $P$ and $<,>$ denotes the canonical scalar product. 

\begin{defi} We define a polyhedral complex $\mathcal{C}(S_\lambda)$ in $\R^n$ as the collection of all cells of the form 
\[ \mathcal{N}_P(F) \cap \{w_0 = -1\}.\]
\end{defi}
Since $\mathcal{N}_P$ is a complete fan, this polyhedral complex has support $\R^n$. If we only consider the cells of the form
\[ \mathcal{N}_P(F) \cap \{w_0 = -1\} \mbox{ for } \dim(F) \geq 1,\]
we get a pure polyhedral complex with support $\mathcal{T}(S_\lambda)$.

\subsection{Compactifications of buildings} 

In general there are different ways to compactify a Bruhat-Tits building $\mathfrak{B}(G,K)$ associated to a reductive group $G$ over $K$. In \cite{bose}, Borel and Serre attach the so-called Tits building at infinity in order to derive finiteness results for arithmetic group cohomology. In \cite{we2} a concrete compactification of the Bruhat-Tits building of $PGL$ is studied, which can be identified with the space of isometry classes of non-Archimedean norms on $K^n$. In \cite{we1} we investigate a dual approach which fits into the description of the building with lattices in $K^n$. 

In \cite{we3} we develop a general theory of compactifications. For every semisimple group $G$ over a non-Archimedean local field $K$ and every faithful, geometrically irreducible algebraic representation $\rho: G \rightarrow GL(W)$ we define a compactification of $\mathfrak{B}(G,K)$. Its boundary can be identified with the union of Bruhat-Tits buildings associated to certain types of parabolics in $G$. The strategy is the following: Use the combinatorics of the weights given by the representation $\rho$ in order to compactify one apartment. Then generalize Bruhat-Tits theory to define groups $P_x$ for all $x$ in the compactified apartment and glue all compactified apartments together as in the definition of $\mathfrak{B}(G,K)$ (see 1.2). 
There is an even more general approach: In \cite{rtw1} we realize the Bruhat-Tits building $\mathfrak{B}(G,K)$ inside the Berkovich analytic space $G^{an}$ and use the projection to analytical flag varieties of $G$ to obtain a family of compactifications of $\mathfrak{B}(G,K)$. This fits together with the approach in \cite{we3}, see \cite{rtw2}.

In this section we  look at the finite family of compactifications of $\mathfrak{B}(SL(V),K)$ defined in \cite{we3}.
Note that the assumption in \cite{we3} that the ground field $K$ is a local field is not needed in the specific examples we discuss here. 

Let $\rho: SL(V) \rightarrow GL(W)$ be a geometrically irreducible algebraic representation of the group $SL(V)$ on a $K$-vector space $W$. For every basis $\Delta$  of the root system $\Phi(T, SL(V))$ we denote by $\mu_0(\Delta)$ the corresponding highest weight of $\rho$. 
\begin{defi}
We define a fan $\mathcal{F}_\rho$ as the set of all faces of the cones
\[C_\Delta(\rho) = \{ x \in A: \mu_0(\Delta) (x) \geq \mu(x) \mbox{ for all weights }\mu \mbox{ of }\rho\},\] 
where $\Delta$ runs over the bases of $\Phi(T, SL(V))$. 
\end{defi}

{\bf Example 1.} Suppose that $\rho = id$. The weights of the identical representation are $\{a_1, \ldots, a_n\}$. For the basis $\Delta =\{ a_{12}, a_{23}, \ldots, a_{n-1 n}\}$ the highest weight is $\mu_0(\Delta) = a_1$. Hence 
\[C_\Delta(\rho)=\{x \in A: a_1(x) \geq a_i(x) \mbox{ for all }i\}= \{ \sum_{i=1}^n x_i \eta_i \in A: x_1 \geq x_i \mbox{ for all }i\}.\] 
Since the Weyl group (which is isomorphic to the symmetric group ${\mathcal{S}}_n$) acts simply transitively on the set of bases, every maximal cone is of the form
\[\Gamma_k = \{\sum_{i=1}^n x_i \eta_i\in A: x_k \geq x_i \mbox{ for all }i\}\]
for some $k = 1,\ldots, n$. Hence the fan $\mathcal{F}_\rho$ consists of the cones $\Gamma_k$ and of all their faces.

{\bf Example 2. } Suppose that $\rho$ is a representation with highest weight $\mu_0(\Delta) = n a_1 + (n-1) a_2 + \ldots + 2 a_{n-1} + a_n$. Then for every root $a_{i i+1}$ there exists a weight $\mu$ such that $\mu_0(\Delta) - \mu = a_{i i+1}$. Hence for $\Delta = \{a_{12}, \ldots, a_{n-1 n}\}$ we find
\[C_\Delta(\rho) = \{ x \in A: a_{12} (x) \geq 0, \ldots, a_{n-1 n}(x) \geq 0 \} = \{\sum_{i=1}^n x_i \eta_i: x_1 \geq x_2 \geq \ldots \geq x_n\}.\]
In this case we write  $C_\Delta(\rho)= \mathfrak{C}(\Delta)$. It is the Weyl cone associated to $\Delta$, see \cite{bou}, chapter VI, \S 1.5.  The fan $\mathcal{F}_\rho$ is the Weyl fan consisting of all $\mathfrak{C}(\Delta)$ and their faces.

Note that for every $\rho$ the cone $C_\Delta(\rho)$ contains the Weyl cone $\mathfrak{C}(\Delta)$. Since the union of all Weyl cones is the total space $A$, we deduce that every fan $\mathcal{F}_\rho$ has support $A$. The Weyl cones for different bases are different. In general however, we have
$C_\Delta(\rho) = C_{\Delta'}(\rho)$, whenever $\mu_0(\Delta) = \mu_0(\Delta')$, and this may happen for  $\Delta \neq \Delta'$.

 The fan $\mathcal{F}_\rho$ can be used to define a compactification $\overline{A}_\rho$ of $A$, see
\cite{we3}, section 2. In fact, we put $\overline{A}_\rho = \bigcup_{C \in \mathcal{F}_\rho} A / \langle C \rangle$ and we endow this space with a topology given by tubular neighbourhoods around boundary points. For a more streamlined definition of this topology see \cite{rtw1}, appendix B.  

Note that the compactification associated to the fan in example 1 above is the one studied in \cite{we1}, whereas the compactification to the fan in example 2 is Landvogt's polyhedral compactification studied in \cite{la}.

The fan $\mathcal{F}_\rho$ and hence the compactification $\overline{A}_\rho$ only depend on the Weyl chamber face containing the highest weight of $\rho$, see \cite{we3}. Hence we obtain a finite family of compactifications of $A$ in this way. Identifying $A$ with $\mathbb{T}^{n-1}$ we also obtain a finite family of compactifications for the tropical torus. The compactification of $\mathbb{T}^{n-1}$ given by example 1 could be regarded as a tropical analog of projective space, see also \cite{jo}. 

Fix a geometrically irreducible representation $\rho$ and let $\lambda$ be a partition such that $\rho$ is isomorphic to $\mathfrak{S}_\lambda(V)$.
We will now compare $\mathcal{F}_\rho$ with the polyhedral complex $\mathcal{C}(S_\lambda)$ given by the tropical Schur polynomial $S_\lambda$, cf. section 2.

\begin{thm}
Assume that all the coefficients of the Schur polynomial $S_\lambda$ have valuation zero. Then the fan 
$\mathcal{F}_\rho$ coincides with the image of $\mathcal{C}(S_\lambda)$ under the map $\R^n \rightarrow \T^{n-1} \simeq A$. 
\end{thm}

{\bf Proof: }Note that our hypothesis implies that all cells in $\mathcal{C}(S_\lambda)$ are invariant under $(x_1,\ldots, x_n) \mapsto (x_1+a, \ldots, x_n+a)$ for any $a \in \R$. It suffices to show that the images of the cones of maximal dimension in $\mathcal{C}(S_\lambda)$  are precisely the cones of maximal dimension in $\mathcal{F}_\rho$, i.e. that they are of the form
$C_\Delta = \{ x \in A: \mu_0(\Delta) (x) \geq \mu(x) \mbox{ for all weights }\mu\}$ for some basis  $\Delta$ of $\Phi(T, SL(V))$.

The cones of maximal dimension in $\mathcal{C}(S_\lambda)$ are by definition of the form $\mathcal{N}_P(v) \cap \{w_0 = -1\}$, where $v$ is a vertex of $P$. Since by assumption all coefficients of $S_\lambda$ have valuation zero, the polygon $P$ can be identified with the convex hull of all $(\nu_1, \ldots, \nu_n)$ with $c_N \neq 0$ in $\R^n$. Fix a  vertex $v$ of this polygon. Then $v = (\mu_1, \ldots, \mu_n) $ with $M = (\mu_1, \ldots, \mu_n)$ such that $c_M \neq 0$, and the corresponding cone in $\mathcal{C}(S_\lambda)$ is equal to
\[\mathcal{C}(v) =\{ (x_1, \ldots, x_n): \sum_i \mu_i x_i \geq \sum_i \nu_i x_i \mbox{ for all }N= (\nu_1, \ldots, \nu_n) \mbox{ with }c_N \neq 0\}.\]
Since this is a cone of full dimension in $\R^n$, its topological interior is non-empty. Hence there exists a point $x \in \R^n$ such that $\sum_i \mu_i x_i > \sum_i \nu_i x_i$  for all $N = (\nu_1, \ldots, \nu_n)$ with $c_N \neq 0$ and $N \neq M$. The image of $x$ in $A$ lies in one of the Weyl cones $\mathfrak{C}(\Delta)$, where $\Delta$ is a basis of the root system. 
Now the set  of coefficients $c_N$ of $S_\lambda$ with $N = (\nu_1, \ldots, \nu_n)$ and $c_N \neq 0$ corresponds bijectively to the set of weights, if we map $c_N$ to $\nu = a_1 \nu_1 + \ldots + a_n \nu_n $. Let $\mu$ be the weight $\mu_1 a_1 + \ldots + \mu_n a_n$. Then $\mu_0(\Delta) - \mu$ is a linear combination of elements in $\Delta$ with non-negative coefficients. Since all elements in $\Delta$ are non-negative on the image of  $x$, we find $\mu_0(\Delta) = \mu$. 
Therefore the image of $\mathcal{C}(v)$ in $A$ is equal to $\mathcal{C}_\Delta$. 
Since every cone in $\mathcal{F}_\rho$ is a face of some $\mathcal{C}_\Delta$, our claim follows.\hfill$\Box$

Let us have a look at the hypothesis that all Kostka numbers occuring as coefficients of $S_\lambda$ have valuation zero. 
Note that this assumption is of course fulfilled if the residue field of $K$ has characteristic zero. Besides, for fixed $\lambda$, it is fulfilled for fields of almost all residue characteristics. 

The theorem implies that the polyhedral complex given by $\mathcal{T}(S_\lambda)$ is in fact a fan,
i.e. all faces are convex cones.  This can also be seen directly from the description of the polyhedral complex in this special case.

\subsection{Tropical convexity}

We consider the situation of example 1 in section 2.2. The identical representation of $SL(V)$ corresponds to the partition $\lambda= (1,0, \ldots, 0)$. Hence by theorem 2.4 the corresponding compactification of the apartment $A$ is given by the polyhedral complex associated to the Schur polynomial $S_{(1,0,\ldots, 0)} (z_1, \ldots, z_n) = z_1 + \ldots+ z_n$.
Transferred to  $\mathbb{T}^{n-1}$  it consists of all faces of the maximal cones 
\[\Gamma_k = \{(x_1, \ldots, x_n) + \R(1, \ldots, 1)  \in \mathbb{T}^{n-1}: x_k \geq x_i \mbox{ for all }i = 1, \ldots, n\}.\]
We have seen that this polyhedral complex arises in the compactification of the building $\mathfrak{B}(SL(V),K)$ associated to the identical representation of $SL(V)$. It also  plays a role in the definition of tropical convexity. 

In \cite{dest}, Develin and Sturmfels define tropical polytopes as tropical convex hulls of finitely many points in $\mathbb{T}^{n-1}$.  Every tropical polytope is the support of a polyhedral complex, see \cite{dest}, section 3, which can be described as follows (see also \cite{jo}).

Let $M=\{v_1,\ldots, v_r\}$ be a finite subset of points in $\T^{n-1}$. 
For every $x \in \mathbb{T}^{n-1}$ define the type of $x$ (depending on $M$) as
\[\mbox{type}(x) = (T_1,\ldots, T_n),\]
where 
\[T_k = \{i: v_i -x \in \Gamma'_k\}\]
with
\[\Gamma'_k = \{ (x_1, \ldots, x_n) + \R(1, \ldots, 1)  \in \mathbb{T}^{n-1}: x_k \leq x_i \mbox{ for all }i = 1, \ldots, n\}.\]  Then every type $T = (T_1, \ldots, T_n)$ defines a polyhedron
\[X_T = \{ x \in \mathbb{T}^{n-1}: x \mbox{ has type }(T_1, \ldots, T_n)\}.\]
The union of all bounded polyhedra of this form is precisely the tropical convex hull of the finite set $M$. 

The cones $\Gamma_k'$ used in tropical convexity are the min-convex twins of the cones $\Gamma_k$. To be precise, the map
$\R \rightarrow \R$ given by $x \mapsto -x$ is an isomorphism between the tropical $(\min, +)$-semiring and the tropical $(\max,+)$-semiring, and we have $-\Gamma'_k = \Gamma_k$. Therefore a $(\max,+)$-version of tropical convexity would use  $\Gamma_k$ instead of $\Gamma'_k$. 

\subsection{Summary and open questions}

We have seen that some features of the Bruhat-Tits building associated to $SL(V)$ can be interpreted in terms of tropical geometry. Every apartment in $\mathfrak{B}(SL(V),K)$ is a tropical projective space. We see the stabilizers of points in the building as tropical stabilizers via a tropical matrix action. Besides, we can identify the fans used to compactify apartments with the polyhedral complex associated to tropical Schur polynomials. It is a natural question if this can be generalized to other reductive groups.

Besides, the compactification of $\mathfrak{B}(SL(V),K)$ described in example 1 in section 2.2 is based on a fan which is also used in tropical convexity. What happens if we keep the group $SL(V)$, but use the fans associated to different compactifications, i.e. the fans associated to  different tropical Schur polynomials? We can then define a version of convexity adapting the definitions from \cite{dest}. What kind of geometry does this lead to?


\small 
\begin{thebibliography}{xxxxxxx}
\bibitem[Bo-Se]{bose} A. Borel, J.-P. Serre: 
{\it Cohomologie d'immeubles et de groupes $S$-arithm\'etiques}. 
Topology { 15} (1976) 211-232.



\bibitem[Br-Ti1]{brti1}F. Bruhat, J. Tits: {\it Groupes r\'eductifs sur 
un corps local. I. Donn\'ees radicielles valu\'ees.}
Publ. Math. IHES { 41} (1972), 5-252.

\bibitem[Br-Ti2]{brti2}F. Bruhat, J. Tits: {\it Groupes r\'eductifs sur 
un corps local. II. Sch\'emas en groups. Existence d'une donn\'e radicielle valu\'ee.} 
Publ. Math. IHES { 60} (1984), 5-184.


\bibitem[Br-Ti3]{brti3} F. Bruhat, J. Tits: {\it Sch\'emas en groupes et immeubles des 
groupes classiques sur un corps local.} Bull. Soc. math. France {112} (1984)
259-301.

\bibitem[Bou]{bou} N. Bourbaki: Lie groups and Lie algebras. Chapter 4-6. Springer 2002. 

\bibitem[De-St]{dest} M. Develin, B. Sturmfels: {\it Tropical convexity}. Doc. Math. { 9} (2004) 1-27. 

\bibitem[EKL]{ekl} M. Einsiedler, M. Kapranow, D. Lind: {\it Non-archimedean amoebas and tropical varieties.} J. Reine Angew. Math. 601 (2006), 139–157.




\bibitem[Gr]{gr} J.A. Green: Polynomial representations of $GL_n$. Lecture Notes in Mathematics 830.  Springer 1980. 

\bibitem[Jos]{jo} M. Joswig: {\it Tropical convex hull computations.} Preprint 2008. ArXiv 0809.4694.

\bibitem[JSY] {jsy} M. Joswig, B. Sturmfels, J. Yu: {\it Affine buildings and tropical convexity}. Albanian J. Math. 1 (2007) 187–211.

\bibitem[La]{la}E. Landvogt: A compactification of the Bruhat-Tits building. 
Lecture Notes in Mathematics {1619}. Springer 1996.

\bibitem[Mac]{mac} I.G. Macdonald: Symmetric functions and Hall polynomials. Oxford Science Publications 1995.

\bibitem[Pro]{pro} C. Procesi: Lie Groups. An Approach through Invariants and Representations. Springer 2007.

\bibitem[RTW1]{rtw1} B. R\'emy, A. Thuillier, A. Werner: {\it Bruhat-Tits theory from Berkovich's point of view - I. Realizations and compactifications of Bruhat-Tits buildings}. Preprint 2009. ArXiv 0903.1245.

\bibitem[RTW2]{rtw2}  B. R\'emy, A. Thuillier, A. Werner: {\it Bruhat-Tits theory from Berkovich's point of view -II. Satake compactifications of buildings}. In preparation. 

\bibitem[Rou]{rou} G. Rousseau: Euclidean buildings. In: G\'eom\'etries \`a courbure n\'egative or nulle, groupes discrets et rigidit\'es (A. Parreau, L. Bessi\`eres and B. R\'emy, eds.) S\'eminaires et Congr\`es no. 18, Soci\'et\'e math\'ematique de France 2008. 


\bibitem[We1]{we1} A. Werner: {\it Compactification of the Bruhat-Tits building of 
PGL by lattices of smaller rank.} Documenta Math. {6} (2001) 315-342.

\bibitem[We2]{we2} A. Werner: {\it Compactification of the Bruhat-Tits building of PGL by seminorms}.
Math. Z. {248} (2004) 511-526.

\bibitem[We3]{we3} A. Werner: {\it Compactifications of {B}ruhat-{T}its buildings associated to linear representations}.
Proc. Lond. Math. Soc. {95} (2007) 497-518. 

\bibitem[Zie]{zi} G. Ziegler: Lectures on Polytopes. Springer 2007. 

\end{thebibliography}
\end{document}